\documentclass[12pt]{amsart}
\usepackage{latexsym}
\usepackage{a4}
\usepackage{amssymb}

\newtheorem{thm}{Theorem}[section] 
 
\newtheorem{lem}[thm]{Lemma} 
 
\newtheorem{cor}[thm]{Corollary} 
\theoremstyle{definition}

\numberwithin{equation}{section} 
\newcommand{\ab}[1]{{\mathbf{#1}}} 
\newcommand{\ob}[1]{{\mathbb{#1}}}

\newcommand{\VecTwo}[2]{ 
   \left( 
   \begin{smallmatrix} 
      #1 \\ #2 
   \end{smallmatrix} 
   \right) 
   } 

\newcommand{\N}{\Bbb{ N}}

\newcommand{\setsuchthat}{\,\, \pmb{|} \,\,}

\newcommand{\vb}[1]{\mathbf{#1}}

\newcommand{\Pol}{\mathsf{Pol}}

\newcommand{\Con}{\mathrm{Con}}

\newcommand{\meet}{\wedge}
\newcommand{\join}{\vee}

\newcommand{\Int}[2]{{I[{#1},{#2}]}}

\newcommand{\OInt}[2]{{\ob{I}[{#1},{#2}]}}

\newcommand{\algop}[2]{( {#1}, {#2} )}

\newcommand{\R}{\mathcal{R}}

\newcommand{\pot}{P}
\newcommand{\epsi}{\varepsilon}
\newcommand{\tup}[3]{(#1_{#2},\dots,#1_{#3})}
\title{Sequences of commutator operations}
\author{Erhard Aichinger}

\address{Erhard Aichinger,
Institut f\"ur Algebra,
Johannes Kepler Universit\"at Linz,
4040 Linz,
Austria}
\email{\tt erhard@algebra.uni-linz.ac.at}

\author{Neboj\v{s}a Mudrinski} 
\address{Neboj\v{s}a Mudrinski,
Department of Mathematics and Informatics,
Faculty of Sciences,
University of Novi Sad,
21000 Novi Sad,
Serbia and
Institut f\"ur Algebra,
Johannes Kepler Universit\"at Linz,
4040 Linz,
Austria}
\email{\tt nmudrinski@dmi.uns.ac.rs}

\subjclass[2000]{06B10 (06A07, 08A40)}
\urladdr{http://www.jku.at/algebra}
\thanks{Supported by the Austrian Science Fund (FWF):P24077}
\keywords{lattices, sequences of operations, commutators}
\date{\today}

\begin{document}
\bibliographystyle{amsalpha}
\begin{abstract}
   Given the congruence lattice $\ob{L}$ of a finite algebra
   $\ab{A}$
   with a Mal'cev term, we look for those sequences 
   of operations on $\ob{L}$ that are sequences of 
   higher commutator operations of expansions of $\ab{A}$.
   
   The properties of higher commutators proved
   so far  delimit the number of such sequences: 
   the number is always at most countably infinite; if it is 
   infinite, then $\ob{L}$ is the union of two
   proper subintervals with
   nonempty intersection.
\end{abstract} 

\maketitle
\section{Introduction}
 It is well known that for a finite algebra with a Mal'cev term,
 the isomorphism type of the congruence lattice 
 yields some information on the binary commutator operation.
 For example, it is well-known that the diamond $\ob{M}_3$ as a congruence
 lattice forces an algebra $\ab{A}$ with Mal'cev term to be abelian, 
 and hence the commutator operation
 to satisfy $[1,1]_{\ab{A}} = 0$. As a consequence of the results of this note,
 the congruence lattice of a finite non-nilpotent Mal'cev algebra
 is equal to the set-theoretic union of two of its proper subintervals; hence
 congruence lattices that are no such union force the algebra to be nilpotent.
 This result is obtained by investigating the higher commutator operations
 as defined by \cite{Bu:OTNO}.
 Given a lattice $\ob{L}$, we will try to delimit
 the number of  sequences $([.,.], [.,.,.],\ldots)$ of
  operations on $\ob{L}$ that could be the sequence of higher commutator operations
  of some Mal'cev algebra with  congruence lattice isomorphic to $\ob{L}$.
   We will see that the number of such sequences is at most countable,
  and we characterise when this number is finite.

This is motivated by the search of a classification of finite algebras with a Mal'cev term modulo
polynomial equivalence. 
We call two universal algebras \emph{polynomially equivalent} if they
are defined on the same universe and they have the same clone of polynomial
operations. For example, for a set $M$ and its power set $\pot(M)$,
the Boolean algebra $\algop{\pot(M)}{\cap,\cup,\mbox{}'}$ and
the Boolean ring $\algop{\pot (M)}{\Delta, \cap}$ are polynomially equivalent
since the fundamental operations of each of these two algebras are 
polynomial operations of the other one.
There are several invariants of an algebra
that depend on the clone of polynomial operations,
but not on the particular choice of fundamental operations.
One of these invariants is the congruence lattice, expanded
with the binary commutator operation introduced in 
\cite{Sm:MV}, cf. \cite{FM:CTFC,MMT:ALVV}. 
Generalizing the binary commutator operation,
A.\ Bulatov introduced multi-placed
 commutators
for an algebra $\ab{A}$ \cite[Definition~3]{Bu:OTNO}.
 For each $k \in \N$, and each $k$-tuple 
$(\alpha_1,\ldots, \alpha_k) \in (\Con (\ab{A}))^k$, he defined a congruence
$[\alpha_1,\ldots,\alpha_k]_{\ab{A}}$ of $\ab{A}$ and named it 
the 
\emph{$k$-ary commutator} of $\alpha_1,\ldots, \alpha_k$.
When $\ab{A}$ has a Mal'cev term, \cite{Mu:Diss,AM:SAOH} discuss
several properties of these higher commutators.
As for the binary term condition commutator,
 these commutator operations are completely
determined by the clone of polynomial functions of an algebra.
So with each algebra $\ab{A}$, we can associate the \emph{commutator
structure} of $\ab{A}$. This is the structure
$\algop{\Con (\ab{A})}{\wedge, \vee,(f_i)_{i \in \N}}$, where 
$f_i : (\Con (\ab{A}))^i \to \Con (\ab{A}), (\alpha_1,\ldots,\alpha_i) 
\mapsto [\alpha_1,\ldots,\alpha_i]_{\ab{A}}$.
If $\ab{A}$ has a Mal'cev term, then $\algop{\Con (\ab{A})}{\wedge, \vee}$
is a modular lattice, and for all $n,k \in \N$
  with $k \le n$, we have
  \begin{itemize}
       \item (HC1) $f_n (\alpha_1,\ldots,\alpha_n) \le \alpha_k$.
       \item (HC2) $f_n$ is monotonous.
       \item (HC3) $f_{n+1} (\alpha_1,\ldots, \alpha_{n+1}) \le f_{n} (\alpha_2, \ldots, \alpha_{n+1})$.
       \item (HC4) $f_n$ is symmetric.
       \item (HC7) $f_n$ is join distributive.
       \item (HC8) $f_{k}(\alpha_{1},\ldots,\alpha_{k-1},f_{n-k+1}\tup{\alpha}{k}{n})\leq f_n\tup{\alpha}{1}{n}$.
  \end{itemize}
These properties have been stated and established in \cite{Bu:OTNO, Mu:Diss,AM:SAOH}.
We note that the properties (HC5) and (HC6) listed in \cite{AM:SAOH} are missing from the list, since
they are not properties of the sequence $(f_i)_{i \in \N}$ but rather on the relation 
of the higher commutators with the underlying algebra.
We call two algebras \emph{commutator equivalent} if they have
the same commutator structure. Since an algebra
has its congruence relations and
its higher commutator operations completely determined by
the  clone of polynomial functions,
polynomially equivalent algebras are commutator equivalent.
For a converse, it is open whether two finite
Mal'cev algebras $\ab{A}$ and $\ab{B}$ 
with $\Pol_3 (\ab{A}) = \Pol_3 (\ab{B})$
and the same commutator structure must be polynomially
equivalent.

Let us now consider an arbitrary sequence $(f_i)_{i \in \N}$ of operations on a lattice $\ob{L}$
such that for each $i \in \N$, the function $f_i$ is an $i$-ary operation on $\ob{L}$.
We say the sequence $(f_i)_{i \in \N}$ is  
\emph{admissible} if it satisfies the list of properties
given above.
In the present note we will investigate the following problem:
 \begin{quote}
     Given a finite modular lattice $\ob{L}$, how many
     admissible sequences can be defined on $\ob{L}$?
\end{quote}
Hence given the isomorphism type of the congruence lattice of a Mal'cev algebra, we want
to delimit the number of possible higher commutator operations
on this algebra.

\section{The cardinality of the set of admissible sequences} 
   Let $\ob{L}$ be a complete lattice, and let $\delta, \epsi \in \ob{L}$. We say that
   $(\delta, \epsi)$ is a \emph{splitting pair} of $\ob{L}$ if 
   $\delta < 1$, $\epsi > 0$, and for all $\alpha \in \ob{L}$, we have
   $\alpha \ge \epsi$ or $\alpha \le \delta$. 
   A splitting pair is called \emph{strong} if $\delta \ge \epsi$.
   The lattice $\ob{L}$ \emph{splits}  if it has a 
   splitting pair, it \emph{splits strongly} if it has
   a strong splitting pair.

   Let us remark that
   this splitting property has often arisen in the following context:
   A splitting pair $(\delta, \epsi)$ in the congruence lattice of 
   an algebra $\ab{A}$ is a rich source of congruence preserving operations on 
   $\ab{A}$ because every finitary operation that is constant
   on $\delta$-classes and maps into one $\epsi$-class is a congruence
   preserving function. \cite{Ai:2ACA, HMP:CTOG} are just two examples
   in which the splitting property of the congruence lattice was used in this way.

      Let  $n \in \N$, and let 
      $f : \ob{L}^n \to \ob{L}$. Then $f$ is \emph{join distributive}
      if for all $i \in \{1,\ldots, n\}$, and for all
      $\alpha_1,\ldots,\alpha_n \in \ob{L}$ and all families $\langle
       \beta_j \setsuchthat j \in J \rangle$ of elements of $\ob{L}$,
       we have
      $f (\alpha_1,\ldots, \alpha_{i-1}, \bigvee_{j \in J} \beta_j, \alpha_{i+1},\ldots, \alpha_n) =
       \bigvee_{j \in J}
       f (\alpha_1,\ldots, \alpha_{i-1}, \beta_j, \alpha_{i+1},\ldots, \alpha_n)$.
      The function $f$ is \emph{symmetric} if
      $f (\alpha_1,\ldots, \alpha_n) = f(\alpha_{\pi(1)}, \ldots, \alpha_{\pi(n)})$ for
      all $\alpha_1,\ldots, \alpha_n \in \ob{L}$ and $\pi \in S_n$, and $f$ is
    \emph{monotonous} if it preserves $\le$.
      The sequence
      $(f_i)_{i \in \N}$ is an \emph{operation sequence on} $\ob{L}$,
      if for all $i \in \N$, $f_i : \ob{L}^i \to \ob{L}$. 
            The operation sequence $(f_i)_{i \in \N}$ is called \emph{admissible}
      if it satisfies the properties (HC1), (HC2),  (HC3), (HC4), (HC7), (HC8).
    
   \begin{thm} \label{thm:cardinality}
        Let $\ob{L}$ be a finite modular lattice. Then
        the number of admissible operation sequences on $\ob{L}$
        is finite if and only if $\ob{L}$ does not split strongly,
        and countably infinite otherwise.
    \end{thm}
    The proof will be completed at the end of Section~\ref{sec:prfcard}.

          Let
      $(f_i)_{i \in \N}$ and $(g_i)_{i \in  \N}$ be operation sequences on the lattice $\ob{L}$.
      We say $(f_i)_{i \in \N} \sqsubseteq (g_i)_{i \in \N}$ if for all $i \in \N$ and for
      all $\alpha_1,\ldots, \alpha_i \in \ob{L}$, we have
      $f_i (\alpha_1,\ldots, \alpha_i) \le g_i (\alpha_1,\ldots, \alpha_i)$.
   
  \begin{thm} \label{thm:infinite}
     Let $\ob{L}$ be a finite lattice, and let $S$ be the set of
     all admissible operation sequences on $\ob{L}$.
     Then $S$ is at most countable, and 
     $(S, \sqsubseteq)$ has no infinite descending chains and no infinite antichains.
  \end{thm}
This result will be proved in Section~\ref{sec:prfcard}.
 
\section{Preliminaries on lattices and ordered sets}

By $\ob{B}_2$, we denote the two element lattice on the set $\{0,1\}$, and by $\ob{M}_3$, we 
denote the diamond.
The lattice $\ob{M}_3$ does not split.
It is easy to see that the lattices $\ob{B}_2$ and $\ob{M}_2 := \ob{B}_2 \times \ob{B}_2$ split, but do not split strongly.
The three element chain $\{0,\theta, 1\}$  splits strongly with splitting pair $(\theta, \theta)$.

\begin{lem} \label{lem:nosplitLat}
   Let $\ob{L}$ be a modular 
   lattice of finite height that does not split strongly.
   Then there are $n \in \N_0$ and a lattice $\ob{M}$ such that 
   $\ob{M}$ does
   not split and
   $\ob{L}$ is isomorphic to $\ob{M} \times (\ob{B}_2)^n$.
\end{lem}
\emph{Proof:}
   We proceed by induction on the height of $\ob{L}$.
   If the height is $0$, then $|\ob{L}| = 1$ and $\ob{L} \cong \ob{L} \times (\ob{B}_2)^0$.
   Now assume that $|\ob{L}| > 1$. If $\ob{L}$ does not 
   split, we take $\ob{M} := \ob{L}$ and $n := 0$.
   Now assume that $\ob{L}$ has a splitting pair 
   $(\alpha, \beta)$. We choose an atom $\epsi$ and 
   a coatom $\delta$ of $\ob{L}$ with
   $\delta \ge \alpha$ and $\beta \ge \epsi$. Then
   $(\delta, \epsi)$ is a splitting pair of $\ob{L}$,
   and since $\ob{L}$ does not split strongly,
   we have $\delta \not\ge \epsi$.
    Let $\ob{L}_1$ be the interval $\OInt{0}{\delta}$,
   and let $\ob{L}_2 := \OInt{0}{\epsi}$.
   By a theorem of Birkhoff \cite[Theorem~2.31]{MMT:ALVV},
   the sublattice of $\ob{L}$ generated by $L_1 \cup L_2$
   is isomorphic to $\ob{L}_1 \times \ob{L}_2$. But since $(\delta, \epsi)$ is a splitting 
   pair, we have $(x \meet \delta) \join (x \meet \epsi) = x$ for
   all $x \in \ob{L}$. To see this, assume $x \le \delta$.
   Then  $(x \meet \delta) \join (x \meet \epsi) = 
           x \join (x \meet \epsi) = x$.
   If $x \ge \epsi$, then $(x \meet \delta) \join (x \meet \epsi) =
                           (x \meet \delta) \join \epsi =
                           x \meet (\delta \join \epsi) = x \meet 1 = x$.
   Hence the sublattice generated by $L_1 \cup L_2$ is equal to $\ob{L}$.
   Thus $\ob{L}$ is isomorphic to $\ob{L}_1 \times \ob{L}_2$.
   
       The lattice $\ob{L}_2$ is isomorphic to $\ob{B}_2$.      
      The lattice $\ob{L}_1$ does not split strongly: suppose $(\delta_1, \epsi_1)$ is a strong splitting
      pair of $\ob{L}_1$. Then $((\delta_1, \epsi), (\epsi_1, 0))$ is a strong splitting pair of 
       $\ob{L}_1 \times \ob{L}_2$, and therefore, 
       $\ob{L}$ has a strong splitting pair, a contradiction.
      Hence applying the induction hypothesis to $\ob{L}_1$, we obtain a lattice
      $\ob{M}$ that does not split and $n \in \N_0$ with $\ob{L}_1 \cong \ob{M} \times \ob{B}_2^{n}$,
      and therefore
      $\ob{L} \cong \ob{L}_1 \times \ob{B}_2 \cong \ob{M} \times \ob{B}_2^{n+1}$. \qed

Let $\ob{A} = (A, \le)$ be a partially ordered set.
We say that $\ob{A}$ satisfies the \emph{descending chain condition}
 if there is no infinitely descending chain
$a_1 > a_2 > a_3 > \ldots$. The \emph{ascending chain condition}
 is defined dually. 
For $m \in \N$, we define a partially ordered
set $\ob{A}^m = (A^m, \le)$, where
$(a_1,\ldots, a_m) \le (b_1,\ldots, b_m)$ if for all $i \in \{1, \ldots, m\}$,
we have $a_i \le b_i$.
For $\ob{A} := \algop{\N}{\le}$, the following lemma is known
as Dickson's Lemma \cite{Di:FOTO}. 
\begin{lem}[{\cite[Lemma~1.2]{La:WASO}, \cite[p.195, Example~(4)]{AH:FTIS}}]  \label{lem:dic}
   Let $\ob{A}$ be a partially ordered set with
   the descending chain condition and no infinite antichains.
   Then $\ob{A}^m$ satisfies the descending chain condition
    and has no infinite antichains.
\end{lem}
   A subset $I$ of $\N_0^m$ is an \emph{upward closed set}
   if for all $\vb{a} \in I$ and $\vb{b} \in \N_0^m$ 
   with $\vb{a} \le \vb{b}$, we have $\vb{b} \in I$.
   It is easy to see that every upward closed set $U$ is
   uniquely determined by its minimal elements. 
   Since the set of minimal elements of $U$ is an antichain,
   Lemma~\ref{lem:dic} implies that this set is finite.
   This establishes the following lemma.
\begin{lem} \label{lem:countable}
   Let $m \in \N$.    There are exactly countably infinitely many upward closed
    subsets of $\N_0^m$.
\end{lem}
We will also use the following theorem from order theory:
\begin{thm} [cf.{\cite[Corollary~4.3]{AH:FTIS},\cite[Theorem~1.2]{Ma:AOMI}}] \label{thm:ucfac}
  Let $m \in \N$, and
  let $\mathcal{L}$ be the set of upward closed
  subsets of $\N_0^m$. Then the partially ordered
  set $(\mathcal{L}, \subseteq)$ has no infinite 
  antichain and no infinite ascending chain.
\end{thm}

\section{Sequences of operations} \label{sec:prfcard}
 First, we prove that the set of admissible operation sequences on a finite lattice
 is at most countable and satisfies certain order properties.
   
\emph{Proof of Theorem~\ref{thm:infinite}:}
    Let $m := |\ob{L}|$,  let $\{\alpha_1,\ldots,\alpha_m\}$ be the
    set of elements of $\ob{L}$, and
    let $F := (f_i)_{i \in \N}$ be an admissible sequence.
    Then for $(a_1,\ldots,a_m) \in \N_0^m \setminus \{ (0,\ldots, 0) \}$, we define
    $E( F, (a_1,\ldots, a_m))$ by
    \[
        E( F, (a_1,\ldots, a_m)) :=   f_{j} (
             \underbrace{\alpha_1, \ldots, \alpha_1}_{a_1 \text{ times}},
              \ldots, 
             \underbrace{\alpha_m, \ldots, \alpha_m}_{a_m \text{ times}} 
                               ),
    \]
    where $j := \sum_{k=1}^m a_k$.

       For $\alpha \in \ob{L}$, we define
   $\R_F (\alpha)$ as the subset of $\N_0^m$ given by
   \begin{multline}
      \R_F (\alpha) =
      \{ (a_1,\ldots,a_m) \in \N_0^m \setminus \{(0, \ldots, 0)\} 
           \setsuchthat
            E (F, (a_1,\ldots, a_m))
          \le \alpha \}.
   \end{multline}
 Since $F$ is an admissible sequence, $\R_F (\alpha)$ is an upward
 closed subset of $\N_0^m$.
 Let $F = (f_i)_{i \in \N}$ and $G= (g_i)_{i \in \N}$ be two admissible
 sequences on $\ob{L}$. We will now show
 that $F \sqsubseteq G$ if and only if for all $\alpha \in \ob{L}$, we have
 $\R_G (\alpha) \subseteq \R_F (\alpha)$.
    For the ``only if''-direction, we let
   $\alpha \in \ob{L}$ and $\vb{a} = (a_1, \ldots, a_m) \in \N_0^m$ such
   that $\vb{a} \in \R_G (\alpha)$.
   Then  
   $E(G, (a_1,\ldots, a_m)) \le \alpha$, and thus, since $F \sqsubseteq G$,
   $E(F, (a_1,\ldots, a_m)) \le  \alpha$, which implies $(a_1,\ldots, a_m) 
   \in \R_F (\alpha)$.
   For the ``if''-direction, we let
   $k \in \N$ and $\beta_1, \ldots, \beta_k \in \ob{L}$.
   Using the symmetry of $f_k$ and $g_k$, we obtain $(a_1, \ldots, a_m) \in \N_0^m$
   such that $f_k (\beta_1, \ldots, \beta_k) =
             E(F, (a_1,\ldots, a_m))$
    and
    $g_k (\beta_1, \ldots, \beta_k) =
            E (G, (a_1,\ldots, a_m))$.
   From the last equality, we obtain that $\vb{a}$ lies 
   in $\R_G (g_k (\beta_1, \ldots, \beta_k))$. Hence
   we have $\vb{a} \in \R_F (g_k (\beta_1, \ldots, \beta_k))$.
   Using the symmetry of $f_k$, this implies
   $f_k (\beta_1, \ldots, \beta_k) \le g_k (\beta_1, \ldots, \beta_k)$. 
   Denoting by
   $\mathcal{U}$ be the set of upward closed subsets
   of $\N_0^m$, we have just proved that
     $(S, \sqsubseteq)$ is isomorphic
   to a subset of 
   the dual of $(\mathcal{U}, \subseteq)^m$.
   Now from  Lemma~\ref{lem:countable}, we obtain that
   $S$ is at most countable.
   By Theorem~\ref{thm:ucfac}, $(\mathcal{U}, \subseteq)$
    has no infinite antichain and no infinite ascending chain.
   Applying Lemma~\ref{lem:dic} to the dual
   of $(\mathcal{U}, \subseteq)$, we obtain that
   $(\mathcal{U}, \subseteq)^m$ satisfies
   the ascending chain condition and has no infinite antichains. Hence
   $(S, \sqsubseteq)$ satisfies the descending chain condition and has no infinite
    antichains.
\qed

  \begin{lem} \label{lem:prod}
    Let $\ob{L}_1, \ob{L}_2$ be 
    lattices, let 
    $\ob{L} := \ob{L}_1 \times \ob{L}_2$, and let
   $(f_i)_{i \in \N}$ be an admissible operation sequence on $\ob{L}$.
    Then for all $n \in \N$, 
    \begin{equation} \label{eq:fj}
        f_n (\VecTwo{\alpha_1}{\beta_1},\ldots, \VecTwo{\alpha_n}{\beta_n}) =
        f_n (\VecTwo{\alpha_1}{0},\ldots, \VecTwo{\alpha_n}{0}) \join 
        f_n (\VecTwo{0}{\beta_1},\ldots, \VecTwo{0}{\beta_n}).
    \end{equation}
  \end{lem}
   \emph{Proof:} We compute 
        $f_n (\VecTwo{\alpha_1}{\beta_1},\ldots, \VecTwo{\alpha_n}{\beta_n}) =
         f_n (\VecTwo{\alpha_1}{0} \join \VecTwo{0}{\beta_1},
              \ldots,
              \VecTwo{\alpha_n}{0} \join \VecTwo{0}{\beta_n})$.
        Using join distributivity, the last expression is equal
        to the join of $2^n$ expressions of the
        form $f (a_1,\ldots, a_n)$ with $a_i \in \{ \VecTwo{\alpha_i}{0},
                                                    \VecTwo{0}{\beta_i} \}$.
         If $a_i = \VecTwo{\alpha_i}{0}$ and $a_j =\VecTwo{0}{\beta_i}$,
         then by (HC1), $f (a_1,\ldots, a_n) \le a_i \meet a_j = 0$.
         Hence $f_n (\VecTwo{\alpha_1}{\beta_1},\ldots, \VecTwo{\alpha_n}{\beta_n})$
         is equal to the join of the two summands in the right hand side of~\eqref{eq:fj}
      \qed
   \begin{lem} \label{lem:01}
       Let $\ob{B}_2$ be the $2$-element lattice.
        Then there are exactly three admissible operation sequences on $\ob{B}_2$;
         these are $(f_i)_{i \in \N}$,
        $(g_i)_{i \in \N}$, and  $(h_i)_{i \in \N}$ with $f_n (\alpha_1,\ldots, \alpha_n) = 0$ for
        all $n \in \N$,
        $g_1 (1) = 1$ and $g_i = f_i$ for $i \ge 2$, 
        and $h_n (\alpha_1,\ldots, \alpha_n) = \alpha_1 \meet \ldots \meet \alpha_n$ for all $n \in \N$.
    \end{lem}
       \emph{Proof:} Let $(s_i)_{i \in \N}$ be an admissible operation sequence
        on $\ob{B}_2$.
         By (HC1), we  have $s_i (\alpha_1,\ldots, \alpha_n) = 0$
         if $0 \in \{\alpha_1, \ldots, \alpha_n \}$. Hence we only need to determine
          $s_i (1,\ldots, 1)$.
         In the case $s_1 (1) = 0$, the property (HC3) yields $s_n (1,\ldots,1) \le
          s_{n-1}(1,\ldots, 1) \le \cdots \le s_1 (1) = 0$ for all $n \in \N$, and
          thus $(s_i)_{i \in \N} = (f_i)_{i \in \N}$.
          In the case that $s_1 (1) = 1$ and $s_2 (1,1) = 0$, we have
          $s_n (1\ldots, 1) = 0$ by repeated application of (HC3). 
           In the case $s_1 (1) = s_2 (1,1) = 1$, (HC8) yields
          $s_n (1,\ldots, 1) \ge s_{n-1} (1,\ldots, 1, s_2(1,1)) 
            = s_{n-1} (1,\ldots, 1)$ for all $n \ge 3$, and thus
          $(s_i)_{i \in \N} = (h_i)_{i \in \N}$.
         \qed
     \begin{lem} \label{lem:nosplit}
        Let $\ob{L}$ be a finite lattice that does not split, let
        $n$ be the number of atoms of $\ob{L}$, and
        let $(f_i)_{i \in \N}$ be an admissible 
         operation sequence on $\ob{L}$. Then for all 
        $k \ge n$, we have $f_k(\gamma_1,\ldots, \gamma_k) = 0$
        for all $\gamma_1,\ldots, \gamma_k \in \ob{L}$.
     \end{lem}
     \emph{Proof:}
      Let $\alpha_1,\ldots,\alpha_n$ be all the atoms of $\ob{L}$. If $n=1$, 
      then $\ob{L}$ splits. 
      Therefore, $n\geq 2$. For each $i \in \{1,\ldots, n\}$, we define
$
A(i):=\{x\in \ob{L}\,|\,x\ngeq\alpha_i\}.
$
We first show that  for all $i \in \{1,2,\ldots, n\}$,
we have $\bigvee \{\alpha \setsuchthat \alpha \in A(i) \}=1$.
Let $\delta := \bigvee \{\alpha \setsuchthat \alpha \in A(i) \}$. Then 
for every $x \in L$, we have $x \ge \alpha_i$ or $x \le \delta$.
Hence if $\delta < 1$, then $(\delta, \alpha_i)$ is a splitting pair,
contradiction the assumptions.
Now if $k \ge n$, using (HC3) and (HC7), we obtain
\begin{multline*}
f_k (1,\dots,1) \le f_n (1,\ldots, 1) 
= f_n (\bigvee_{x_1 \in A(1)} x_1, \ldots, \bigvee_{x_n \in A(n)} x_n)   \\
=\bigvee_{\tup{x}{1}{n}\in A(1)\times\dots\times A(n)}f_n \tup{x}{1}{n}.
\end{multline*}
We will now show that each $f_n \tup{x}{1}{n}$ is equal to $0$.
Suppose $f_n \tup{x}{1}{n} > 0$. Then there is an atom $\alpha_j \in L$
with $f_n \tup{x}{1}{n} \ge \alpha_j$.
Hence $\alpha_j \le f_n \tup{x}{1}{n} \le x_j$. This contradicts
$x_j \in A(j)$. \qed

This Lemma has a consequence on the congruence lattice of a non-nilpotent algebra.
An  algebra $\ab{A}$ with a Mal'cev term is \emph{nilpotent} if and only
if its lower central series of congruences defined by $\gamma_1 := 1$, $\gamma_{n}  := [1, \gamma_{n-1}]_{\ab{A}}$
for $n \ge 2$ reaches $0$ after finitely many steps.
We recall that a direct product $\ab{B} = \ab{A}_1 \times \cdots \times \ab{A}_n$ is 
\emph{skew-free} if for every congruence relation $\beta$ of $\ab{B}$, there
are congruences $\alpha_1 \in \Con (\ab{A}_1), \ldots,
                 \alpha_n \in \Con (\ab{A}_n)$ such that
for all $(a_1,\ldots, a_n), (b_1,\ldots, b_n) \in B$, we have
$((a_1,\ldots, a_n), (b_1,\ldots, b_n)) \in \beta$ if and only 
if $(a_i, b_i) \in \alpha_i$ for all $i \in \{1, \ldots, n\}$.
\begin{cor}
   Let $\ab{A}$ be a finite algebra with a Mal'cev term.
   Then we have:
   \begin{enumerate}
      \item \label{it:sn1}
         If $\ab{A}$ is not nilpotent, 
         then its congruence lattice $\Con (\ab{A})$ splits.
      \item \label{it:sn2}
          If $\Con (\ab{A})$ does not split strongly,
   then there exist $n \in \N_0$ and algebras $\ab{B}, \ab{C}_1, 
   \ldots, \ab{C}_n$ such that
    $\ab{A}$ is isomorphic to the direct product 
    $\ab{B} \times \ab{C}_1 \times \cdots \times \ab{C}_n$,
    $\ab{B}$ is nilpotent,  each $\ab{C}_i$ is simple,
    and the direct product is skew-free.
   \end{enumerate}
\end{cor}  
\emph{Proof:}
    \eqref{it:sn1}
    Assume that the lattice $\Con (\ab{A})$ does not split. 
    Then by Lemma~\ref{lem:nosplit}, there is an $n \in \N$ such 
    that the $n$-ary higher commutator operation of $\ab{A}$
    satisfies $[1,\ldots,1]_{\ab{A}} = 0$.
    By (HC8) and (HC2), we obtain that then the $n$-th term
    $\gamma_n$ of the lower central series of $\ab{A}$ 
    satisfies $\gamma_n = 0$. Hence $\ab{A}$ is nilpotent,
    contradicting the assumptions.
  
    \eqref{it:sn2}
    We  assume that the congruence lattice of $\ab{A}$ 
    does not split strongly.
    Then Lemma~\ref{lem:nosplitLat}
    yields an $n \in \N_0$ and a lattice $\ob{M}$ that does not
    split such that
    $\Con (\ab{A})$ is isomorphic via some isomorphism
    $\iota$ to $\ob{M} \times \ob{B}_2^n$.
    For $i \in \{0,\ldots, n\}$, let
     $\nu_i := \iota^{-1}  ( (1, 1, \ldots, 1, 0 , 1, \ldots, 1))$ with
     $0$ at the $(i+1)$-th place. Using the fact that these congruences
     permute, we obtain (cf. \cite[p.161]{MMT:ALVV}) that
     $\ab{A}$ is isomorphic to $\prod_{i=0}^n (\ab{A} / \nu_i)$.
     Since $\Con (\ab{A} / \nu_0)$ is isomorphic to  
     $\ob{M}$, the congruence lattice of $\ab{A}/ \nu_0$ does not split,
      and hence, by the the first part of this corollary,
      $\ab{A}/\nu_0$ is nilpotent. For $i \ge 1$,
      $\nu_i$ is a coatom of $\Con (\ab{A})$ and  
       $\ab{A}/ \nu_i$ is simple.
     Hence $\ab{B} := \ab{A}/\nu_0$ and $\ab{C}_i := \ab{A}/\nu_i$
     satisfy $\ab{A} \cong \ab{B}  \times \prod_{i=1}^n \ab{C}_i$.
     For every $\theta \in \Con (\ab{A})$, we have
     $\theta = \iota^{-1} (\iota (\theta)) 
             = \iota^{-1} ((\iota (\theta) \join (0,1,1 \ldots, 1)) 
                           \meet
                            \ldots
                           \meet
                           (\iota(\theta) \join (1,1,1,\ldots, 0)))
             = (\theta \join \nu_0) \meet \ldots \meet (\theta \join \nu_n)$,
    and therefore the direct product is skew-free by \cite[Lemma~IV.11.6]{BS:ACIU}.
   \qed

\begin{thm} \label{thm:finite} 
     Let $\ob{L}$ be a finite modular lattice, and let $S$ be the
     set of all admissible sequences on $\ob{L}$. Then $S$ is infinite if
     and only if $\ob{L}$ splits strongly.
\end{thm}
\emph{Proof:}
       Let us assume that $\ob{L}$ does not split strongly.
       Then by Lemma~\ref{lem:nosplitLat}, $\ob{L}$ is 
       isomorphic to a direct product $\ob{M} \times \ob{B}_2^n$
       such that $\ob{M}$ does not split.
       Now by Lemmas~\ref{lem:nosplit}~and~\ref{lem:01}, 
       on each of the direct factors, there are only
       finitely many admissible operation sequences, and thus
       by Lemma~\ref{lem:prod}, $S$ is finite.
      
       If $\ob{L}$ splits strongly, then we choose
        a strong splitting pair $(\delta, \epsi)$, and we define
       an operation sequence $(f_i)_{i \in \N}$ by $f_1 (\alpha_1) := \alpha_1$
       for all $\alpha_1 \in \ob{L}$, and for $i \ge 2$,
       $f_i (\alpha_1, \ldots, \alpha_i) := 0$ if there
       exists an $j \in \{1,\ldots, i\}$ with $\alpha_j \le \delta$,
       and $f_i (\alpha_1, \ldots, \alpha_i) := \epsi$  else.
        Let $g_i (\alpha_1, \ldots, \alpha_i) = 0$ for $i \in \N$.
        Now we show that for each $k \in \N$, the sequence
        $(h^{(k)}_i)_{i \in \N}$ defined by
        $h^{(j)}_i := f_i$ for $i \le j$ and
        $h^{(j)}_i := g_i$ for $ i > j$ is an admissible sequence.
       
        We first show that each $f_i$ satisfies (HC1).
        Supposing that (HC1) fails for some $\alpha_1, \ldots, \alpha_i$, 
        we have
        $f_i (\alpha_1, \ldots, \alpha_i) = \epsi$ and thus
       $\alpha_j \not\le \delta$ for all $j \in \{1, \ldots, i\}$.
        Thus $\alpha_j \ge \epsi$ for all $j$, and therefore
         $f_i (\alpha_1, \ldots, \alpha_i) \le \bigwedge_{j=1}^i \alpha_j$.
        (HC2), (HC3), and (HC4) are immediate consequences of the definitions.
         Now for join distributivity, having already established
        (HC1--4), we only need to prove
        \[
           f_i( \bigvee_{j \in J} \beta_j, \alpha_2, \ldots, \alpha_i)
           \le \bigvee_{j \in J}  f_i (\beta_j, \alpha_2, \ldots, \alpha_i)
        \]
         for all families $\langle \beta_j \setsuchthat j \in J \rangle$
         from $\ob{L}$.
         Suppose that the right hand side is $0$.
         Then either one of the $\alpha_k$ satisfies
         $\alpha_k \le \delta$, implying that the left hand
         side is $0$, or all $\alpha_k$ satisfy
         $\alpha_k \not\le \delta$. Then we have
         $\beta_j \le \delta$ for all $j \in J$. This implies
         $\bigvee_{j \in J} \beta_j \le \delta$, and therefore the
         left hand side is $0$ as well.

         In order to prove (HC8) for each sequence
         $(h^{(j)})_{i \in \N}$, we observe that for
         $j \le i - 2$, and for
         every nested expression
          of the form $f_{j+1} (\alpha_1,\ldots, \alpha_j, f_{j-i} (\alpha_{j+1}, \ldots,
                            \alpha_i))$, we have                          
          $f_{j+1} (\alpha_1,\ldots, \alpha_j, f_{j-i} (\alpha_{j+1}, \ldots,
                            \alpha_i)) \le
           f_{j+1} (\alpha_1, \ldots, \alpha_j, \epsi) = 0$.
        \qed                         

Now Theorem~\ref{thm:cardinality} follows immediately from Theorems~\ref{thm:infinite}~and~\ref{thm:finite}.
As a consequence, we give an upper bound on the number of pairwise
commutator inequivalent algebras with Mal'cev term on a finite
universe.
\begin{cor}
    Let $A$ be a finite set, let $I$  be an infinite set, let $\ob{L}$
    be  a sublattice
    of the lattice of equivalence relations on $A$, and let 
    $(\ab{B}_i)_{i \in I}$ be  a family of algebras with universe $A$
     such that for  each  $i \in I$, $\ab{B}_i$ has a Mal'cev term $d_i$ and
     $\Con (\ab{B}_i) = \ob{L}$,
     and for all $i, j \in I$ with $i \neq j$,  
     $\ab{B}_i$ and $\ab{B}_j$ are not
     commutator equivalent.
     Then $|I| \le \aleph_0$, and $\ob{L}$
     is the union of two intervals $\Int{0}{\delta} \cup \Int{\epsi}{1}$
     with $0 < \epsi \le \delta < 1$.
\end{cor}
\emph{Proof:}
     For each $i \in I$ and $j \in \N$, we define
     $h^{(i)}_j (\alpha_1,\ldots,\alpha_j) := 
      [\alpha_1, \ldots, \alpha_j]_{\ab{B}_i}$.
     Since each $\ab{B}_i$ has a Mal'cev term, each
     $(h^{(i)}_j)_{j \in \N}$ is an admissible sequence.
     Since all $\ab{B}_i$ are commutator inequivalent,
     we get an infinite set of admissible sequences. Thus by Theorem~\ref{thm:cardinality},
     $I$ is countably infinite and $\ob{L}$ splits strongly. \qed  

\providecommand{\bysame}{\leavevmode\hbox to3em{\hrulefill}\thinspace}
\providecommand{\MR}{\relax\ifhmode\unskip\space\fi MR }
\providecommand{\MRhref}[2]{%
  \href{http://www.ams.org/mathscinet-getitem?mr=#1}{#2}
}
\providecommand{\href}[2]{#2}

\end{document}